\documentstyle[12pt,a4]{article}

\def\rixrel#1#2{\mathrel{\mathop{#2}\limits_{#1}}}

\def\p{\sf p}

\def\m{\sf m}

\def\P{\bf P}
\def\Z{\bf Z}

\begin{document}

\title{Bounds on cohomology and Castelnuovo-Mumford regularity}
\vspace{10mm}
\author{Chikashi Miyazaki\thanks{This author would like to thank Massey
University for
financial support and the Department of Mathematics for its friendly
atmosphere while
writing this paper.}\\Nagano National College of Technology\\716 Tokuma,
Nagano 381\\
Japan\\e-mail: miyazaki@ei.nagano-nct.ac.jp\\\
\and and\\\ \\Wolfgang Vogel\\
Department of Mathematics\\Massey University\\Palmerston North\\New
 Zealand
\\ e-mail: W.\@Vogel@massey.ac.nz}

\date{}
\maketitle
\newpage

\section{Introduction}

Let $X\subseteq{\P}_K^N$ be a projective scheme over an algebraically
closed field $K$ . We denote by ${\cal I}_X$ the ideal sheaf of
 $X$. Then
$X$ is said to be
$m$-regular if $H^i({\P}_K^N,{\cal I}_X(m-i))=0$ for all $i\geq
 1$ (cf.
\cite{Mum}). The
Castelnuovo-Mumford regularity $\mbox{reg}(X)$ of $X\subseteq
{\P}_K^N$,
firstly introduced by
Mumford by generalizing ideas of Castelnuovo, is the least such integer
$m$. The interest in this
concept stems partly from the well-known fact that $X$ is $m$-regular
 if
and only if for every
$p\geq 0$ the minimal generators of the $p$-th syzygy module of the
defining ideal $I$ of
$X\subseteq{\P}_K^N$ occur in $\mbox{degree}\leq m+p$ (see, e.g.,
\cite{EG}). There are
good bounds in some cases if $X$ is assumed to be smooth (see \cite{BEL}
and the references there).

Our interest is to consider the case $X$ being locally Cohen-Macaulay and
equidimensional.
Under this assumption there
is a non-negative integer $k$ such that
$$(X_0,\cdots,X_N)^k\left[\rixrel{\ell\in{\Z}}{\oplus}H^i\left(
{\P}_K^N,{\cal I}_X(\ell)\right)\right]=0
\mbox{\ for\ }1\leq i\leq\dim X,$$
where ${\P}_K^N=\mbox{Proj\ } K[X_0,\cdots,X_N]$. In this case
$X\subseteq{\P}_K^N$ is called a
$k$-Buchsbaum scheme. A refined version is a $(k,r)$-Buchsbaum scheme
introduced in \cite{FV},
\cite{HMV}, \cite{HV}:

Let $k$ and $r$ be integers with $k\geq 0$ and $1\leq r\leq \dim
 X$. Then
we call $X\subseteq{\P}_K^N$
a $(k,r)$-Buchsbaum scheme if, for all $j=0,\cdots,r-1$, $X\cap
 V$ is a
$k$-Buchsbaum scheme for every
$(N-j)$-dimensional complete intersection $V$ in ${\P}_K^N$ with
$\dim(X\cap V)=\dim(X)-j$,
(see (2.3) and see also (2.4) for an equivalent definition in case
 $k=1$.)

In recent years upper bounds on the Castelnuovo-Mumford regularity
 of such
a variety
$X\subseteq {\P}_K^N$ have been given by several authors in terms
 of
$\dim(X)$, $\deg(X)$, $k$
and $r$, (see, e.g., \cite{HMV}, \cite{HM}, \cite{HV}, \cite{NS1},
\cite{NS2}). These bounds
are stated as follows:
$$\mbox{reg}(X)\leq\left\lceil\frac{\deg(X)-1}{\mbox{codim\
}(X)}\right\rceil+C(k,r,d),$$
where $d=\dim X$, $C(k,r,d)$ is a constant depending on $k$, $r$
 and $d$,
and $\left\lceil n\right\rceil$
is the smallest integer $\ell\geq n$ for a rational number $n$. In
 case $X$
is arithmetically
Cohen-Macaulay, that is, $k=0$, it is well-known that $C(k,r,d)\leq
 1$,
(see, e.g., \cite{SV2},
\cite{SV3}). We assume $k\geq 1$. In case $r=1$ it was shown that
$C(k,1,d)\leq\left(\begin{array}{c}
d+1\\2\end{array}\right)k-d+1$ in \cite{HV}. In \cite{HM} (in a slightly
weaker form) and
\cite{NS1} it was improved to $C(k,1,d)\leq(2d-1)k-d+1$. Further,
 a better bound
$C(k,1,d)\leq dk$ was obtained in \cite{NS2}. The general case $r\geq
 1$
was firstly studied in
\cite{HMV}, which was improved in \cite{HV} by showing that
$C(k,r,d)\leq(r-1)k+\left(\begin{array}{c}
d+2-r\\2\end{array}\right)k-d+1$.

The purpose of this paper is to give bounds on $\mbox{reg\ }(X)$
 in terms
of $\dim (X)$, $\deg (X)$,
$k$ and $r$, which improve some of the previous results. In general case
 we show that
$$C(k,r,d)\leq dk-r+1$$
in (3.2). Moreover, in case $k=1$, we show that
$$C(1,r,d)\leq\left\lceil\frac{d}{r}\right\rceil.$$
in (3.3)

Our methods here is to use a spectral sequence theory for graded
 modules
developed in \cite{M1},
\cite{M2}, \cite{M3} in order to get bounds on the local cohomology
 and the
Castelnuovo-Mumford
regularity.

Let $R$ be the coordinate ring of $X\subseteq{\P}_K^N$ with the
 homogeneous
maximal ideal $\m$. Then we define $a_i(R)=\sup\{n;H_{\m}^i(R)_n\neq 0\}$
for
$i=0,\cdots,d+1$, and $\mbox{reg\ }(R)=\max\{a_i(R)+i;i=0,\cdots,d+1\}
$. It
is easy to see that
$\mbox{reg\ }(X)=\mbox{reg\ }(R)+1$. Among $a_i(R)$'s, $a_{d+1}(R)$,
 which
is called an
$a$-invariant of $R$ (cf. \cite{GW}), plays an important role to
characterize the graded ring $R$.
The results (2.7) and (2.9) on bounds on the cohomology in terms
 of
$a$-invariants are applied to
(3.2) and (3.3).

In section 2, we consider bounds on the cohomology for graded
 modules.
First we explain a
spectral sequence theory in order to get cohomological bounds. Theorem
(2.5) is the technical key
results of this paper. This theorem gives several important results
 stated
in (2.7), (2.8), (2.9),
(2.10), (2.11).

In section 3, we describe bounds on the Castelnuovo-Mumford regularity
 for
locally
Cohen-Macaulay varieties $X\subseteq{\P}_K^N$. Theorem (3.2)
 and
Theorem (3.3)
are main results of this section.

In section 4, we construct sharp examples of (2.8) and (2.11), and
 conclude by
describing some open problems.

\section{Bounds on local cohomology}

Throughout this paper let $R=K[R_1]$ be a Noetherian graded $K$-algebra,
where $K$ is an infinite
field. We denote by $\m$ the maximal homogeneous ideal of $R$.
 Let $M$
be a graded $R$-module.
We write $[M]_n$ for the $n$-th graded piece of $M$ and $M(p)$ for
 the
graded module with
$[M(p)]_n=[M]_{p+n}$. We set:
$$\sup(M)=\sup\{n;[M]_n\neq 0\}.$$
If $M=0$, we set $\sup(M)=-\infty$. Now we assume that $M$ is
 a finitely
generated graded
$R$-module with $\dim M=d>0$. Then we set:
$$a_i(M)=\sup(H^i_{\m}(M)),\qquad i=0,\cdots,d.$$
Also $a_d(M)$ is called the $a$-invariant of $M$ and written as $a(M)$.
 The
Castelnuovo-Mumford
regularity of $M$ is defined as follows:
$$\mbox{reg\ }(M)=\max\{a_i(M)+i; i=0,\cdots,d\}.$$
In order to state our results we use the following notation:
$$\mbox{reg}_n(M)=\max\{a_i(M)+i;i=n,\cdots,d\},$$
where $0\leq n\leq d$. Then we see $\mbox{reg\
}(M)=\mbox{reg}_0(M)=\mbox{reg}_{\mbox{depth\ }(M)}
(M)$.

Let $X$ be a closed subscheme of ${\P}_K^N$. Let $R$ be the coordinate
ring of $X$. Then the
Castelnuovo-Mumford regularity of $X$ is defined as follows:
$$\mbox{reg\ }(X)=\mbox{reg\ }(R)+1.$$

  From now on we assume that the graded module $M$ is a generalized
Cohen-Macaulay module. In other
words $M$ is locally Cohen-Macaulay, that is, $M_{\p}$ is
Cohen-Macaulay for all non-maximal
homogeneous prime ideal $\p$ of $R$, and $M$ is equi-dimensional,
 that is,
$\dim\left(R/{\p}\right)=d$ for all minimal associated prime
 ideals
$\p$ of $M$.
\vspace{5mm}

\noindent{\bf Definition 2.1.}
{\em Let $k$ be a non-negative integer. The graded $R$-module $M$
 is called
a $k$-Buchsbaum module
if ${\m}^kH_{\m}^i(M)=0$ for all $i<d$.}
\vspace{5mm}

\noindent{\bf Definition 2.2.}
(\cite{FV}, \cite{HMV}, \cite{HV}). {\em Let $k$ and $r$ be integers
 with
$k\geq 0$ and $1\leq r\leq d$.
$M$ is called a $(k,r)$-Buchsbaum module if for every s.o.p.
$x_1,\cdots,x_d$ of $M$ we have
$${\m}^kH_{\m}^i(M/(x_1,\cdots,x_j)M)=0,$$
for all non-negative integers $i$, $j$ with $j\leq r-1$ and $i+j<d$.}
\vspace{5mm}

\noindent{\bf Remark.} $(0,d)$-, $(1,d)$-, and $(k,1)$-Buchsbaum
 modules
are the Cohen-Macaulay,
Buchsbaum, and $k$-Buchsbaum modules respectively. The $(1,r)$-Buchsbaum
 module, introduced as
$r$-Buchsbaum module in \cite{M1}, was studied also in \cite{M2},
 \cite{M3}.
\vspace{5mm}

\noindent{\bf Remark.} For a fixed integer $r$, $1\leq r\leq d$,
 every
generalized
Cohen-Macaulay graded module $M$ has $(k,r)$-Buchsbaum property for
 $k$
large enough, see
\cite{HV}, (2.3).
\vspace{5mm}

\noindent{\bf Definition 2.3.}
(\cite{M3}). {\em Let $r$ and $n$ be integers with $1\leq r\leq n\leq
 d$.
Let $x_1,\cdots,x_n$ be a
part of s.o.p. for $M$. We say $x_1,\cdots,x_n$ is $r$-standard,
 if for any
choice
$x_{i_1},\cdots,x_{i_\ell}$ $(\ell\leq r-1)$
$$(x_1,\cdots,x_n)H_{\m}^j(M/(x_{i_1},\cdots,x_{i_\ell})M)=0$$
for $j+\ell<d$.}
\vspace{5mm}

\noindent{\bf Remark.} Standard s.o.p. was introduced in several
 papers,
see, e.g., \cite{T}.
The $r$-standardness is its generalization. A similar generalization
 was
introduced in
\cite{NS1}, (4.2).
\vspace{5mm}

First we note a characterization of $(1,r)$-Buchsbaum modules in
 terms of
s.o.p. of degree one
which improves, e.g., [14], (2.5) of \cite{GM}.
\vspace{5mm}

\renewcommand{\labelenumi}{\theenumi)}
\renewcommand{\theenumi}{\arabic{enumi}}

\noindent{\bf Proposition 2.4.}{\em
Let $R=K[R_1]$ be a graded ring over an infinite field $K$. Let
 $\m$
be the homogeneous maximal
ideal of $R$ generated by $X_0,\cdots,X_N$. Let $M$ be a generalized
Cohen-Macaulay graded
$R$-module with $\dim(M)=d$. We assume that for all integers $0\leq
i_1<\cdots<i_d\leq N$,
$X_{i_1},\cdots,X_{i_d}$ is a s.o.p. for $M$. Let $r$ be an integer
 with
$1\leq r\leq d$. Then the
following are equivalent:
\begin{enumerate}
\item $M$ is a $(1,r)$-Buchsbaum module.
\item For all integers $0\leq i_1<\cdots<i_d\leq N$,
$X_{i_1},\cdots,X_{i_d}$ is an
$r$-standard s.o.p. for $M$.
\end{enumerate}}
\vspace{5mm}

\noindent{\bf Proof.}
It is clear that 1) implies 2). Now we will show that 2) implies
 1) by
induction on $r$.
It is clear in case $r=1$. We assume that $r>1$. To prove that
 $M$ is a
$(1,r)$-Buchsbaum module
we have only to show that, for any part of a s.o.p. $y_1,\cdots,y_r$
 for
$M$, $y_1,\cdots,y_r$ is
$r$-standard, which is equivalent to saying that
$\varphi_{y_1\wedge\cdots\wedge y_r}^q(M)$ is
a zero map for all $r-1\leq q\leq d-1$ by \cite{M3}, (2.3). (Also
 see the
notation in
\cite{M3}.) We set $y_j=\sum_{i=0}^N a_{ji}X_i$ for some $a_{ji}\in
 R$.
Since $X_{i_1},\cdots,
X_{i_r}$ is $r$-standard for all $0\leq i_1<\cdots<i_r\leq N$, we
 see
$\varphi_{X_{i_1}\wedge\cdots\wedge X_{i_r}}^q(M)$ is a zero map by
\cite{M3}, (2.3).
By virtue of \cite{M3}, (3.4), we have that $\varphi_{y_1\wedge\cdots\wedge
y_r}^q(M)$
is a zero map. Thus the assertion is proved.
\vspace{5mm}

The following theorem is the technical key result of this paper.
 The
spectral sequence theory as
developed in \cite{M1}, \cite{M2}, \cite{M3} plays an important role
 to
prove it.
\vspace{5mm}

\noindent{\bf Theorem 2.5.}{\em
Let $M$ be a generalized Cohen-Macaulay graded $R$-module with $\dim
(M)=d>0$. Let $i$ and $n$ be
integers with $0\leq i\leq d-1$ and $1\leq n\leq d$. Let $x_1,\cdots,x_n$
be a part of s.o.p. for
$M$ with $\deg(x_j)=e_j\geq 1$, $j=1,\cdots,n$. We set
$c_j=\max\{e_{i_1}+\cdots+e_{i_j};
1\leq i_1<\cdots<i_j\leq n\}$ for $j=1,\cdots,n$. Then we have:
\begin{enumerate}
\item If $n+i\leq d$, then
$$\begin{array}{lrr}
\lefteqn{\sup(H_{\m}^i(M)/(x_1,\cdots,x_n)H_{\m}^i(M))}\\[2mm]
& & \leq\max\{a_{i+j}(M)+c_{j+1},a_i(M/(x_1,\cdots,x_n)M);j=1,\cdots,n-1\}
{}.
\end{array}$$
\item If $n+i>d$, then
$$\begin{array}{lrr}
\lefteqn{\sup(H_{\m}^i(M)/(x_1,\cdots,x_n)H_{\m}^i(M))}\\[2mm]
& & \leq\max\{a_{i+j}(M)+c_{j+1};j=1,\cdots,d-i\}.
\end{array}$$
\end{enumerate}
Furthermore let $r$ be an integer with $1\leq r\leq
 n$.
Assume that the sequence $x_1,\cdots,x_n$
is
$r$-standard. Then we have:
\begin{enumerate}
\setcounter{enumi}{2}
\item If $n+i\leq d$, then
$$a_i(M)\leq\max\{a_{i+j}(M)+c_{j+1},a_i(M/(x_1,\cdots,x_n)M);
j=r,\cdots,n-1\}.$$
\item If $n+i>d$, then
$$a_i(M)\leq\max\{a_{i+j}(M)+c_{j+1},a(M)+c_{d-i+1};j=r,\cdots,d-i-1\}.$$
\end{enumerate}}
\vspace{5mm}

In order to prove (2.5) we need to apply the following spectral sequence
 corresponding to the graded
$R$-module $M$:

Let $I^\bullet$ be the minimal injective resolution of $M$ in the
 category
of the graded $R$-modules.
Let $K_\bullet$ be the Koszul complex $K_\bullet((x_1,\cdots,x_n);R)$
 for a
part of a s.o.p.
$x_1,\cdots,x_n$ for $M$. Then we consider the double complex
$B^{\bullet\bullet}=\mbox{Hom}_R
(K_\bullet,I^\bullet)$. The filtration $F_t(B^{\bullet\bullet})=\sum_{p\ge
q
t}B^{p,q}$
gives a spectral sequence $\{E_u^{p,q}\}$ (see, e.g., \cite{God}).
 Then we
have the following
isomorphisms:
$$E_1^{p,q}\cong K^p((x_1,\cdots,x_n);H_{\m}^q(M))$$
and
$$H^{p+q}(B^{\bullet\bullet})\cong\left\{\begin{array}{ll}
H^{p+q}((x_1,\cdots,x_n);M), & p+q<n\\[2mm]
H_{\m}^{p+q-n}(M/(x_1,\cdots,x_n)M)(e_1+\cdots+e_n), & p+q\geq
n\end{array}\right.$$
for all $p$, $q$, see, e.g., \cite{M3}. Note that the spectral sequence
$\{E_u^{p,q}\}$ converges
to $H^{p+q}(B^{\bullet\bullet})$. Then we have the following lemma
 by using
the above notation.
\vspace{5mm}

\noindent{\bf Lemma 2.6.}{\em
Let $r$ be an integer with $1\leq r\leq n$. Then the following conditions
are equivalent:
\begin{enumerate}
\item $x_1,\cdots,x_n$ is $r$-standard.
\item $d_s^{p,q}:E_s^{p,q}\rightarrow E_s^{p+s,q-s+1}$ is a zero
 map for
all $p$, $q$, $s$ with
$q\neq d$, $1\leq s\leq r$.
\end{enumerate}}
\vspace{5mm}

\noindent{\bf Proof.}
The lemma follows by induction on $r$ using (3.3) of \cite{M3}.
\vspace{5mm}

Now let us begin to prove (2.5).
\vspace{5mm}

\noindent{\bf Proof of (2.5).}
Let us consider the spectral sequence $\{E_r^{p,q}\}$ corresponding
 to the
graded $R$-module $M$
discussed above for a part of s.o.p. $x_1,\cdots,x_n$ for $M$. We
 see that
$$E_2^{n,i}\cong H_{\m}^i(M)/(x_1,\cdots,x_n)H_{\m}^i(M)(c_n)$$
and
$$H^{n+i}\cong\left\{\begin{array}{ll}
H_{\m}^i(M/(x_1,\cdots,x_n)M)(c_n), & n\leq d-i\\
0, & n>d-i.\end{array}\right.$$
Since $E_1^{p,q}=0$ for $p>n$, we see $E_{\infty}^{p,q}=0$ for
 $p>n$. Thus we
see that the graded
$R$-homomorphism $E_\infty^{n,i}\rightarrow H^{n+i}$ is injective
 and that
$E_\infty^{n,i}$ is a
quotient of the graded $R$-module $E_2^{n,i}$.

First let us prove 1). Take an integer $\ell$ satisfying
$\ell>a_{i+j}(M)+c_{j+1}$ for all
$j=1,\cdots,n-1$ and $\ell>a_i(M/(x_1,\cdots,x_n)M)$. We want to
 show
$\left[H_{\m}^i(M)/
(x_1,\cdots,x_n)\right.$ $\left.H_{\m}^i(M)\right]_\ell=0$.
Since
${[H^{n+i}]}_{\ell-c_n}={[H_{\m}^i
(M/(x_1,\cdots,x_n)M)]}_{\ell}  =  0$, we have
${[E_\infty^{n,i}]}_{\ell-c_n}=0$. On the other hand
we have isomorphisms
\begin{eqnarray*}
{\left[E_1^{n-j-1,i+j}\right]}_{\ell-c_n} & \cong &
{\left[K^{n-j-1}((x_1,\cdots,x_n);
H_{\m}^{i+j}(M))\right]}_{\ell-c_n}\\
& \cong & {\left[\rixrel{1\leq i_1<\cdots<i_{n-j-1}\leq n}{\oplus}
 H_{\m}^{i+j}(M)
(e_{i_1}+\cdots+e_{i_{n-j-1}})\right]}_{\ell-c_n}\\
& \cong &
\rixrel{1\leq k_1<\cdots<k_{j+1}\leq n}{\oplus}{\left[H_{\m}^{i+j}(M)\right]}
_{\ell-(e_{k_1}+\cdots+e_{k_{j+1}})}.
\end{eqnarray*}
Thus we have ${[E_1^{n-j-1,i+j}]}_{\ell-c_n}=0$ for all $j=1,\cdots,n-1$,
 because
$\ell-(e_{k_1}+\cdots+
e_{k_{j+1}})\geq \ell-c_{j+1}>a_{i+j}(M)$. This leads to the equality
${[E_2^{n,i}]}_{\ell-c_n}
={[E_\infty^{n,i}]}_{\ell-c_n}$. Hence we have
${[E_2^{n,i}]}_{\ell-c_n}=0$. Thus the assertion is
proved.

Next let us prove 2). Take an integer $\ell$ satisfying
$\ell>a_{i+j}(M)+c_{j+1}$ for all
$j=1,\cdots,d-i$. Since we have
$${[E_1^{n-j-1,i+j}]}_{\ell-c_n}\cong\rixrel{1\leq k_1<\cdots<k_{j+1}\leq
n}{\oplus}
{[H_{\m}^{i+j}(M)]}_{\ell-(e_{k_1}+\cdots+e_{k_{j+1}})}$$
as shown in the proof of 1), we get ${[E_1^{n-j-1,i+j}]}_{\ell-c_n}=0$
 for
all $j=1,\cdots,d-i$.
Thus we have ${[E_2^{n,i}]}_{\ell-c_n}={[E_\infty^{n,i}]}_{\ell-c_n}=0$.
Hence the assertion
is proved.

To prove 3) we take an integer $\ell$ satisfying $\ell>a_{i+j}(M)+c_{j+1}$
for all $j=r,\cdots,n-1$
and $\ell>a_i(M/(x_1,\cdots,x_n)M)$. Since
${[H^{n+i}]}_{\ell-c_n}={[H_{\m}^i(M/(x_1,\cdots,x_n)
M)]}_\ell=0$, we have ${[E_\infty^{n,i}]}_{\ell-c_n}=0$. Also
 we see
$E_2^{n,i}\cong H_{\m}^i
(M)(c_n)$, because the sequence $x_1,\cdots,x_n$ is 1-standard. Since
$${[E_1^{n-j-1,i+j}]}_{\ell-c_n}\cong \rixrel{1\leq k_1<\cdots<k_{j+1}\leq
n}{\oplus}
{[H_{\m}^{i+j}(M)]}_{\ell-(e_{k_1}+\cdots+e_{k_{j+1}})}.$$
as shown in the proof of 1), we have ${[E_1^{n-j-1,i+j}]}_{\ell-c_n}=0$
 for
all $j=r,\cdots,n-1$.
By using Lemma (2.6), we therefore have
${[E_2^{n,i}]}_{\ell-c_n}={[E_\infty^{n,i}]}_{\ell-c_n}$.
Hence we have ${[H_{\m}^i(M)]}_\ell=0$. Thus the assertion
 is proved.

Finally let us prove 4).  Take an integer $\ell$ satisfying
$\ell>a_{i+j}(M)+c_{j+1}$ for
all $j=r,\cdots,d-i-1$ and $\ell>a(M)+c_{d-i-1}$. Similarly we
 have
${[E_\infty^{n,i}]}_{\ell-c_n}=0$,
$E_2^{n,i}\cong H_{\m}^i(M)(c_n)$, ${[E_1^{n-d+i-1,d}]}_{\ell-c_n}=0$
 and
${[E_1^{n-j-1,i+j}]}_{\ell-c_n}=0$ for all $j=r,\cdots,d-i$.
 By using (2.6)
we therefore
have ${[E_2^{n,i}]}_{\ell-c_n}={[E_\infty^{n,i}]}_{\ell-c_n}$.
 Hence we
have ${[H_{\m}^i(M)]}_\ell
=0$. Thus the assertion is proved.

This completes the proof of (2.5).
\vspace{5mm}

The first application of (2.5) is the following theorem which gives
 bounds
on local
cohomology.
\vspace{5mm}

\noindent{\bf Theorem 2.7.}{\em
Let $M$ be a generalized Cohen-Macaulay graded $R$-module with
$\dim(M)=d>0$. Let
$k_1,\cdots,k_{d-1}$ be non-negative integers satisfying that
${\m}^{k_j}H_{\m}^j(M)=0$
for $j=1,\cdots,d-1$. Let $r$ be an integer with $1\leq r\leq d$.

 Let $x_1,\cdots,x_d$ be a s.o.p. for $M$ with $\deg(x_j)=1$ for
 $j=1,\cdots,d$.
Assume that the sequence $x_1^{\mu_1},\cdots,x_r^{\mu_r}$ is $r$-standard
for some positive
integers $\mu_1,\cdots,\mu_r$. Then we have
\begin{enumerate}
\item $a_i(M)\leq
a_{d-r}(M/(x_1,\cdots,x_r)M)+(d-r-i)+\sum_{\alpha=i}^{d-r-1}k_\alpha+
\sum_{\alpha=1}^r\mu_\alpha-r$ for $i=1,\cdots,d-r-1$.
\item $a_i(M)\leq
a_i(M/(x_1,\cdots,x_{d-i})M)+\sum_{\alpha=1}^{d-i}\mu_\alpha-(d-i)$
 for
$i=d-r,\cdots,d-1$.
\item $\mbox{\em reg}_i(M)\leq
a_{d-r}(M/(x_1,\cdots,x_r)M)+(d-r)+\sum_{\alpha=i}^{d-r-1}k_\alpha
+\sum_{\alpha=1}^r\mu_\alpha-r$ for $i=1,\cdots,d-1$.
\end{enumerate}}
\vspace{5mm}

\noindent{\bf Proof.}
First we will prove 2). By virtue of (2.5.3) we have
$$a_i(M)\leq a_i(M/(x_1^{\mu_1},\cdots,x_{d-i}^{\mu_{d-i}})M)$$
for $i=d-r,\cdots,d-1$, because $x_1^{\mu_1},\cdots,x_{d-i}^{\mu_{d-i}}$
 is
$(d-i)$-standard.
Hence we have
$$a_i(M)\leq
a_i(M/(x_1,\cdots,x_{d-i})M)+\sum_{\alpha=1}^{d-i}\mu_\alpha-(d-i)$$
for $i=d-r,\cdots,d-1$, by using \cite{NS1}, (6.5). Next we will
 prove 1).
Note that
$$a_j(M)\leq\sup(H_{\m}^j(M)/{\m}H_{\m}^j(M))+k_j-1$$
for $j=1,\cdots,d-1$. (See, e.g., \cite{NS2}, (3.3)). By virtue
 of (2.5.2)
we see
\begin{eqnarray*}
a_i(M) & \leq & \sup(H_{\m}^i(M)/(x_1,\cdots,x_d)H_{\m}^i(M))+k_i-1\\[2mm]
& \leq & \max\{a_{i+j}(M)+j+k_i;j=1,\cdots,d-i\}
\end{eqnarray*}
for $i=1,\cdots,d-1$. Now we prove the following claim.
\vspace{5mm}

\noindent{\bf Claim:} Let $i$ be an integer with $1\leq i\leq d-1$.
 Then we have
$$a_i(M)\leq\max\left\{a_j(M)+(j-i)+\sum_{\alpha=i}^{t-1}k_\alpha;j=t,
\cdots
,d\right\}$$
for $t=i+1,\cdots,d$.
\vspace{5mm}

We will show our claim by induction on $t$. The case $t=i+1$ has
 already
shown. In case
$a_t(M)=-\infty$, in other words, $k_t=0$, the claim follows
 immediately
from the hypothesis of induction.
In case $a_t(M)>-\infty$, in other words, $k_t>0$, we have
by the hypothesis of induction
\begin{eqnarray*}
a_i(M) & \leq & \max\left\{a_j(M)+(j-i)+\sum_{\alpha=i}^{t-1}k_\alpha;
j=t,\cdots,d\right\}\\
& \leq &
\max\left\{a_{t+s}(M)+s+k_s+(t-i)+\sum_{\alpha=i}^{t-1}k_\alpha,\right.\\
& & \quad \left.
a_j(M)+(j-i)+\sum_{\alpha=i}^{t-1}k_\alpha;s=1,\cdots,d-t\mbox{\
 and\
}j=t+1,\cdots,d\right\}\\
& = & \max\left\{a_j(M)+(j-i)+\sum_{\alpha=i}^tk_\alpha;j=t+1,\cdots,d
\right\}.
\end{eqnarray*}
Thus the claim is proved.
\vspace{5mm}

In particular, we have
$$a_i(M)\leq\max\left\{a_j(M)+(j-i)+\sum_{\alpha=i}^{d-r-1}k_\alpha;j
=d-r,
\cdots,d\right\}$$
for $i=1,\cdots,d-r-1$. By using 2), we have
\begin{eqnarray*}
a_i(M) & \leq & \max\left\{a_j(M/(x_1,\cdots,x_{d-j})M)+(j-i)+\sum_{\alpha
=i}^{d-r-1}k_\alpha\right.\\
 & & \qquad
\left.+\sum_{\alpha=1}^{d-j}\mu_\alpha-(d-j);j=d-r,\cdots,d\right\}\\
& \leq &
a_{d-r}(M/(x_1,\cdots,x_r)M)+(d-r-i)+\sum_{\alpha=i}^{d-r-1}k_d+\sum_{\alpha
=1}^r
\mu_\alpha-r
\end{eqnarray*}
for $i=1,\cdots,d-r-1$. Hence the assertion 1) is proved. The assertion
 3)
follows immediately
from 1) and 2).

This completes the proof of (2.7).
\vspace{5mm}

Our theorem (2.7) has an important corollary which is used in order
 to
prove (3.2).
\vspace{5mm}

\noindent{\bf Corollary 2.8.}{\em
Let $M$ be a generalized Cohen-Macaulay graded $R$-module with $\dim
(M)=d>0$. Let $k$ and $r$
be integers with $k\geq 1$ and $1\leq r\leq d$. Assume that $M$ is
 a
$(k,r)$-Buchsbaum module.
For $i=1,\cdots,d$,
$$\mbox{\em reg}_i(M)\leq a_{d-r}(M/(x_1,\cdots,x_r)M)+(d-r)+(d-i)k-r$$
for any part of s.o.p. $x_1,\cdots,x_r$ for $M$ with $\deg(x_j)=1$,
$j=1,\cdots,r$.}
\vspace{5mm}

\noindent{\bf Proof.}
Note that $x_1^k,\cdots,x_r^k$ is $r$-standard. By (2.7.3) we have
\begin{eqnarray*}
\mbox{reg}_i(M) & \leq &
a_{d-r}(M/(x_1,\cdots,x_r)M)+(d-r)+\sum_{\alpha=i}^{d-r-1}k+
\sum_{\alpha=1}^r k-r\\
& = & a_{d-r}(M/(x_1,\cdots,x_r)M)+(d-r)+(d-i)k-r.
\end{eqnarray*}
\vspace{5mm}

\noindent{\bf Remark.} The example (4.2) shows that the bound stated
 in (2.8)
is sharp in case $k=r=1$.
\vspace{5mm}

Another application of (2.5) is the following result.
\vspace{5mm}

\noindent{\bf Theorem 2.9.}{\em
Let $M$ be a generalized Cohen-Macaulay graded $R$-module with $\dim
(M)=d>0$. Let $r$ be an
integer with $1\leq r\leq d$. Let $x_1,\cdots,x_d$ be a s.o.p. for
 $M$ with
$\deg(x_j)=1$ for
$j=1,\cdots,d$.
\begin{enumerate}
\item Assume that $x_1,\cdots,x_{d-i}$ is $\min(r,d-i)$-standard
 for a
fixed integer
$i=0,\cdots,d-1$. Then we have $$a_i(M)\leq
a_i(M/(x_1,\cdots,x_{d-i})M)+\left\lceil\frac{d-i}{r}\right\rceil-1.$$
\item If $x_1,\cdots,x_d$ is $r$-standard, then
$$\mbox{\em reg}_i(M)\leq
a_i(M/(x_1,\cdots,x_{d-i})M)+i+\left\lceil\frac{d-i}{r}\right\rceil-1$$
for $i=0,\cdots,d-1$.
\end{enumerate}}
\vspace{5mm}

\noindent{\bf Proof.}
By (2.5.3) we have
\begin{eqnarray*}
a_i(M) & \leq & \max\left\{a_j(M)+(j-i)+1, a_i(M/(x_1,\cdots,x_{d-i})M);
 \right.\\
& & \qquad\left. j=i+r,\cdots,d-1\right\}\\
& \leq & \max\left\{a_j(M)+(j-i)+2,a_i(M/(x_1,\cdots,x_{d-i})M),
 \right.\\
& & \qquad\left.a_{i+r}(M/(x_1,\cdots,x_{d-r-i})M)+r+1;
j=i+2r,\cdots,d-1\right\}\\
& \leq & \max\left\{a_j(M)+(j-i)+2, a_i(M/(x_1,\cdots,x_{d-i})M)+1;
 \right.\\
& & \qquad\left. j=i+2r,\cdots,d-1\right\}.
\end{eqnarray*}
By repeating this step we finally have
\begin{eqnarray*}
a_i(M) & \leq &
\max\left\{a_j(M)+(j-i)+\left\lceil\frac{d-i}{r}\right\rceil-1,\right.\\
& & \qquad
a_i(M/(x_1,\cdots,x_{d-i})M)+\left\lceil\frac{d-i}{r}\right\rceil-2;\\
& & \qquad \left.
j=i+\left(\left\lceil\frac{d-i}{r}\right\rceil-1\right)r,\cdots,d-1\right\}\\
& \leq & a_i(M/(x_1,\cdots,x_{d-i})M)+\left\lceil\frac{d-i}{r}\right\rceil-1
\end{eqnarray*}
for $i=0,\cdots,d-1$. Hence the assertion 1) is proved. The assertion
 2) is
an easy consequence
of 1).
\vspace{5mm}

\noindent{\bf Corollary 2.10.}{\em
Let $M$ be a generalized Cohen-Macaulay graded $R$-module with
$\dim(M)=d>0$. Let $r$ be an integer
with $1\leq r\leq d$. Assume that $M$ is a $(1,r)$-Buchsbaum module.
 For
$i=0,\cdots,d-1$,
$$\mbox{\em reg}_i(M)\leq
a_{d-r}(M/(x_1,\cdots,x_r)M)+(d-r)+\left\lceil\frac{d-i}{r}\right\rceil-1$$
for any part of s.o.p. $x_1,\cdots, x_r$ with $\deg(x_j)=1$, $j=1,\cdots
,r$.}
\vspace{5mm}

\noindent{\bf Proof.} It follows immediately from (2.9).
\vspace{5mm}

\noindent{\bf Remark.} We have assumed that the degree of a s.o.p.
$x_1,\cdots,x_d$ is one in
(2.7), (2.8), (2.9) and (2.10). We took that
assumption for simplicity, although  generalized results are similarly
 proved.
\vspace{5mm}

Further, we have two more results of Theorem (2.5).
\vspace{5mm}

\noindent{\bf Corollary 2.11.}{\em
Let $M$ be a $k$-Buchsbaum graded $R$-module with $\dim(M)=d>0$.
 Then we have
$$\mbox{\em reg}_i(M)\leq a(M)+d+k(d-i)$$
for $i=1,\cdots,d$.}
\vspace{5mm}

\noindent{\bf Proof.}
It follows immediately from the claim in the proof of (2.7).
\vspace{5mm}

\noindent{\bf Remark.} Corollary (2.11) was firstly obtained in \cite{NS2},
(3.4). We will show in this paper that the inequality of (2.11) is
 sharp
for all $d$ and $k$, even in the case $M=R$ is an integral domain,
 see (4.1).
\vspace{5mm}

\noindent{\bf Proposition 2.12}. {\em
Let $M$ be a $(1,r)$-Buchsbaum graded $R$-module with $\dim(M) = d > 0$.
 Then
we have
$$\mbox{\em reg}_i(M)\leq a(M)+d+\left\lceil\frac{d-i}{r}\right\rceil$$
for $i=1,\cdots,d$.}
\vspace{5mm}

\noindent{\bf Proof.} The proof is similar to that of (2.9) by using
 (2.5.4).
\vspace{5mm}

\noindent{\bf Remark.}
Let us consider in (2.12) the special case that $M$ is Buchsbaum,
 that is,
$r=d$. Then
we get for $i\geq 1$
$$\mbox{reg}_i(M)\leq a(M)+d+1.$$
This result can be obtained also by the structure theorem of maximal
Buchsbaum modules. Let us
explain this approach suggested to us by S. Goto: First we take a
polynomial ring $T$
such that $M$ is a maximal $T$-module. By the structure theorem \cite{G},
we see that $M$ as a
graded $T$-module is a direct sum of some twistings of $i$-th syzygy
modules $E_i$. On the
other hand we know that $\mbox{reg\ }(E_i)=i$ and $a(E_i)=i-1$
 for
$i=1,\cdots,d-1$. Hence we
have $$\mbox{reg\ }_i(M)\leq a(M)+d+1$$
for all $i\geq 1$.

\section{Bounds on Castelnuovo-Mumford regularity}

Let $K$ be an algebraically closed field . Let
 $X$ be a
nondegenerate closed subvariety of ${\P}_K^N$ with $\dim(X)=d$.
 Let $R$
be the coordinate ring
of $X$ with $\dim(R)=d+1$. We assume that $X$ is irreducible and
 reduced,
that is, $R$ is an
integral domain.

Before stating our main results of this section we need the following
 lemma
which is well-known,
(see, e.g., \cite{N}, Corollary 2;[20], 4.6(b)).
\vspace{5mm}

\noindent{\bf Lemma 3.1.}{\em
Under the above condition, we have
$$a(R)+d+1\leq \left\lceil\frac{\deg(X)-1}{\mbox{\em codim}(X)}\right\rceil$$}
\vspace{5mm}

The following is our main result of this section, which extends \cite{NS2},
(4.8) and improves
\cite{HMV}, (3.1.5) and \cite{HV}, (3.6).
\vspace{5mm}

\noindent{\bf Theorem 3.2}. {\em
Let $X$ be a non-degenerate closed subvariety of ${\P}_K^N$ with
$\dim(X) = d$ over an algebraically
closed field $K$. Let $R$ be the coordinate
 ring of
$X$. Assume that $X$
is a $(k,r)$-Buchsbaum variety for some integer $k$ and $r$ with
 $k\geq 1$
and $1\leq r\leq d$.
Then we have
$$\mbox{\em reg\ }(X)\leq\left\lceil\frac{\deg (X)-1}{\mbox{\em codim\
}(X)}\right\rceil+(d+1-\mbox{\em depth\ }(R))k
-r+1$$}
\vspace{5mm}

\noindent{\bf Proof.}
It follows immediately from (2.8) and (3.1).
\vspace{5mm}

\noindent{\bf Remark.}
Comparing with the inequality of \cite{HV}, (3.6), we see that our
 result
improves their result
in all cases. We have only to check that
$$dk-r\leq\left((r-1)+\left(\begin{array}{c}d+2-r\\2\end{array}\right)\right
)k-d.$$
Firstly we can easily show that
$$d\leq (r-1)+\left(\begin{array}{c}d+2-r\\2\end{array}\right).$$
So we have only to study the case $k=1$, which is left to the
 readers.
\vspace{5mm}

The next theorem generalizes known results for quasi-Buchsbaum varieties,
see, e.g., \cite{HM},
\cite{NS1}.
\vspace{5mm}

\noindent{\bf Theorem 3.3}. {\em
Let $X$ be a non-degenerate closed subvariety of ${\P}_K^N$ with
$\dim(X) = d$ over an
algebraically closed field $K$. Let $R$ be
 the
coordinate ring of $X$.
Assume that $X$ is a $(1,r)$-Buchsbaum variety for some integer $r$
 with
$1\leq r\leq d$.
Then we have $$\mbox{\em reg\ }(X)\leq\left\lceil\frac{\deg(X)-1}{\mbox{\em
codim\ }(X)}\right\rceil+
\left\lceil\frac{d+1-\mbox{\em depth\ }(R)}{r}\right\rceil.$$}
\vspace{5mm}

\noindent{\bf Proof.}
It follows immediately from (2.10) and (3.1).

\section{Examples and open problems}

The purpose of this section is to describe sharp examples of (2.11)
 and
(2.8). Moreover, we conclude with some open problems.

Let $X\subseteq{\P}_K^N$ be a projective variety with $\dim (X)=d\geq
1$ over an algebraically
closed field $K$ of characteristic $0$. Let $R$ be the coordinate
 ring of
$X$. Assume that $X$ is
$(k,r)$-Buchsbaum for some integers $k$ and $r$ with $k\geq 1$ and
 $1\leq
r\leq d$. Let $L$ be an
$(N-r)$-plane in ${\P}_K^N$ with $\dim (X\cap L)=d-r$. Let
 $R'$ be the
coordinate ring of
$X\cap L$.

Under the above condition, Corollary (2.11), Corollary (2.8) and
Theorem (3.2) are stated as follows:
$$\begin{array}{rl}
\mbox{(2.11)}\qquad & \mbox{reg\ }(X)\leq a(R)+(d+1)+kd+1,\\[5mm]
\mbox{(2.8)}\qquad & \mbox{reg\ }(X)\leq a(R')+(d+1-r)+kd-r+1\\[5mm]
\mbox{(3.2)}\qquad & \mbox{reg\
}(X)\leq\displaystyle\left\lceil\frac{\mbox{deg}(X)-1}{\mbox{codim\
}(X)}\right\rceil+kd-r+1.
\end{array}$$
Note that:
$$a(R)+(d+1)\leq a(R')+(d+1-r)\leq\left\lceil\frac{\deg(X)-1}{\mbox{codim\
}(X)}\right\rceil.$$
Hence (3.2) follows from (2.8). Moreover, (2.11) gives the following bound
in case $r=1$: ${\rm reg}(X) \leq a(R')+d+kd+1$. But this result is
improved in (2.8).

Now we will give sharp examples of (2.11) and (2.8) in Example
(4.1) and Example (4.2) respectively. These examples are based on
 ideas of
\cite{M1}, (3.4) and \cite{M2}, (3.9).

The first example shows that the inequalities of (2.11) are sharp
 for all
$d$ and $k$.
\vspace{5mm}

\noindent{\bf Example 4.1.}
Let $d$ be an integer with $d\geq 1$. Let $Y_j$ be the projective
 line
${\P}_K^1$ over an
algebraically closed field $K$ for $j=1,\cdots,d+1$. Let $Y$ be
 the Segre
product of $Y_j$
($j=1,\cdots,d+1$), that is, $Y=Y_1\times\cdots\times Y_{d+1}$.
 Let
$p_j:Y\rightarrow Y_j$
be the projection for $j=1,\cdots,d+1$. We write an invertible
 sheaf
$p_1^\ast{\cal O}_{{\P}^1}
(n_1)\otimes\cdots\otimes p_{d+1}^\ast{\cal O}_{{\P}^1}(n_{d+1})$
 on
$Y$ as ${\cal O}_Y(n_1,\cdots,
n_{d+1})$ through the isomorphisms $\mbox{Pic\ }(Y)\cong \mbox{Pic\
}(Y_1)\oplus\cdots\oplus
\mbox{Pic\ }(Y_{d+1})\cong{\Z}^{\oplus d+1}$. On the other hand
 $Y$ is
embedded in the projective
space ${\P}_K^N$, where $N=2^{d+1}-1$. Then, for any $n\in{\Z}$,
${\cal O}_{{\P}_K^N}
(n)|_Y\cong{\cal O}_Y(n,\cdots,n)$. There exists an irreducible smooth
effective divisor $X$ of $Y$
corresponding to the invertible sheaf ${\cal O}_Y(n_1,\cdots,n_{d+1})$
 for
all positive integers
$n_1,\cdots,n_{d+1}$. (See, e.g., \cite{Ha}, p231.) Let $k$ be a
 positive
integer. Let us take integers
$n_j=1+(k+1)(j-1)$ for $j=1,\cdots,d+1$. Then we can take a non-singular
subvariety $X$ of $Y$
such that the ideal sheaf ${\cal J}_{X/Y}$ is isomorphic to ${\cal
O}_Y(-n_1,\cdots,-n_{d+1})$.
Let $R$ be the coordinate ring of $X$ in ${\P}_K^N$. Then $R$
 is an
integral domain with
$\dim R=d+1$. In order to get $a(R)$ and $a_i(R)$, $i=1,\cdots,d$,
 we use
the following exact
sequence and the isomorphisms:
$$0\rightarrow H_{\m}^{d+1}(R)\rightarrow \rixrel{\ell\in{\Z}}
{\oplus}H^{d+1}(Y,
{\cal J}_{X/Y}(\ell))\rightarrow\rixrel{\ell\in{\Z}}
{\oplus}H^{d+1}(Y,{\cal O}_Y(\ell))
\rightarrow 0$$
and $$H_{\m}^i(R)\cong\rixrel{\ell\in{\Z}}{\oplus}H^i(Y,{\cal
J}_{X/Y}(\ell)),\quad
i\neq 0,\ d+1,$$
because $Y$ is arithmetically Cohen-Macaulay (see, e.g., \cite{M1},
 (3.1)). Since
we have by K\"{u}nneth's formula:
$$H^{d+1}(Y,{\cal J}_{X/Y}(\ell))\cong H^1({\cal O}_{
{\P}^1}(\ell-n_1))\otimes\cdots\otimes
H^1 ({\cal O}_{{\P}^1}(\ell-n_{d+1}))$$
and $$H^{d+1}(Y,{\cal O}_Y(\ell))\cong H^1({\cal O}_{
{\P}^1}(\ell))^{\otimes d+1},$$
we see that $H^{d+1}(Y,{\cal J}_{X/Y}(\ell))\neq 0$ if $\ell\leq\min_{1\leq
j\leq d+1} (n_j-2)=-1$
and that $H^{d+1}(Y,{\cal O}_Y(\ell))\neq 0$ for $\ell\leq -2$. Thus
 we
have $a(R)=-1$ from
the above exact sequence. Next let us take a fixed integer $i$ with
 $1\leq
i\leq d$. Similarly
we have by K\"{u}nneth's formula
\begin{eqnarray*}
H^i(Y,{\cal J}_{X/Y}(\ell)) & \cong & H^0({\cal O}_{
{\P}^1}(\ell-n_1))\otimes\cdots\otimes
H^0({\cal O}_{{\P}^1}(\ell-n_{d-i+1}))\\
& & \ \ \ \otimes H^1({\cal O}_{
{\P}^1}(\ell-n_{d-i+2}))\otimes\cdots\otimes H^1({\cal O}_{{\P}^1}
(\ell-n_{d+1})).
\end{eqnarray*}
So we see that $H^i(Y,{\cal J}_{X/Y}(\ell))\neq 0$ if
$n_{d-i+1}\leq\ell\leq n_{d-i+2}-2$,
that is, $(k+1)(d-i+1)-k\leq\ell\leq (k+1)(d-i+1)-1$. Thus we have
$a_i(R)=(k+1)(d-i+1)-1$
for $1\leq i\leq d$. Also we have $\mbox{reg}_i(R)=k(d-i+1)+d$
 for $1\leq
i\leq d$, so
$\mbox{reg\ }(R)=(k+1)d$. Further we can easily see that $R$ is
 a
$k$-Buchsbaum ring.
Thus we have a non-singular projective $k$-Buchsbaum variety $X$
 in ${\P}^N$ with
$\dim (X)=d$, $a(R)=-1$ and $\mbox{reg\ }(X)=kd+d+1$. Hence
 this example
yields the
equality stated in (2.11).
\vspace{5mm}

The second example shows that the inequalities of (2.8) are sharp
 in case
$k=r=1$. Also it is possible to give examples even in case $r>1$
 and $k=1$.
In fact we have only
to take some higher dimensional projective space instead of
${\P}_K^1$
in (4.2).
\vspace{5mm}

\noindent{\bf Example 4.2.}
In Example (4.1), we take $k=1$. Then the quasi-Buchsbaum ring
 $R$ gives
the equality of
Corollary (2.8) for all $d$ in case $k=r=1$. In fact, we can
 easily see that
$a(R/hR)=1$ for generic elements $h$ of $R_1$. Then we see that
$$\mbox{reg}_i(R)=2 d-i+1$$
and $$a(R/hR)+d+(d+1-i)\cdot 1-1=2 d-i+1$$
for $i=1,\cdots,d$. Hence this example yields the equality stated
 in (2.8).
\vspace{5mm}

Finally we conclude by describing some open problems which are
 arised
naturally from our
investigation. We take the notation and assumption of the beginning
 of this
section.
\vspace{5mm}

\noindent{\bf Problem 1.}
Generalize (2.11) for $(k,r)$-Buchsbaum varieties for all $r\geq
 1$, or
construct example of $(k,d)$-Buchsbaum varieties with $\dim (X)=d$
satisfying the equality of
(2.11) for all $k$ and $d$.
\vspace{5mm}

We note that the examples in (4.1) are not $(k,d)$-Buchsbaum for
 $d\geq 2$
by using arguments
of \cite{M1}.
\vspace{5mm}

\noindent{\bf Problem 2.}
Improve the bound stated in (2.8) in order to get sharp examples
 for all $k$
and $d$ at least in case $r=1$.

\end{document}